\documentclass[]{article}  
\pdfoutput=1
\usepackage{graphicx}
\usepackage{float}
\usepackage{amsthm}
\usepackage{amssymb}
\usepackage{amsmath}
\usepackage{mathrsfs}
\usepackage{indentfirst}
\usepackage{geometry}
\usepackage{authblk}
\newtheorem{Theorem}{Theorem}[section]
\theoremstyle{definition}

\newtheorem{Lemma}[Theorem]{Lemma}

\theoremstyle{definition}

\title{A complete characterization of graphs for which  $m_G(-1) = n - d - 1 $}
\author{Songnian Xu, Wenhao Zhen, Dein Wong\thanks{Corresponding author. E-mail address: wongdein@163.com.}}
\affil{\textit{School of Mathematics, China University of Mining and Tecnology, Xuzhou, China.}}
\date{}

\begin{document}
\baselineskip 17pt

\title{A complete characterization of graphs for which  $m_G(-1) = n - d - 1 $}

\author{Songnian Xu\\
{\small  Department of Mathematics, China University of Mining and Technology}\\
{\small Xuzhou, 221116, P.R. China}\\
{\small E-mail: xsn1318191@cumt.edu.cn}\\ \\
Wang\thanks{Corresponding author}\\
{\small Department of Mathematics, China University of Mining and Technology}\\
{\small Xuzhou 221116, P.R. China}\\
{\small E-mail: }}

\date{}
\maketitle

\begin{abstract}
Let $G$ be a simple connected graph of order $n$ with diameter $d$.
Let $m_G(-1)$ denote the multiplicity of the eigenvalue $-1$ of the adjacency matrix of $G$, and let $P = P_{d+1}$ be the diameter path of $G$.
If $-1$ is not an eigenvalue of $P$, then by the interlacing theorem, we have $m_G(-1)\leq n - d - 1$.
In this article, we characterize the extremal graphs where equality holds.
Moreover, for the completeness of the results, we also characterize the graphs $G$ that achieve $m_G(-1) = n - d - 1$ when $-1$ is an eigenvalue of $P$.
Thus, we provide a complete characterization of the graphs $G$ for which $m_G(-1) = n - d - 1$.
\end{abstract}

\let\thefootnoteorig\thefootnote
\renewcommand{\thefootnote}{\empty}

\section{Introduction}
In this paper, we consider only simple, connected and finite graphs.
A simple undireted graph $G$ is denoted by $G=(V(G),E(G))$, where $V(G)$ is the vertex set and $E(G)$ is the edge set.
A graph $H$ is called a subgraph of a graph $G$ if $V(H)\subseteq V(G)$ and $E(H)\subseteq E(G)$.
Further, if any two vertices of $V(H)$ are adjacent in $H$ if and only if they are adjacent in $G$, we say that $H$ is an induced subgraph of $G$ and denote this relation as $H\leq G$.
The order of $G$ is the number of vertices of $G$, denoted by $|G|$.
For $V_1\subseteq V(G)$ and $H\leq G$, $N_H(V_1)=\cup_{v\in V_1}\{u\in V(H) | uv\in E(G)\}$.
Let $K\subseteq V(G)$, traditionally, the subgraph of $G$ induced by $K$, written as $G[K]$.
And we sometimes write $G-K$ or $G-G[K]$ to denote $G[V(G)\backslash V(K)]$.
A vertex $v$ of $G$ is said to be pendant if $d(v)=1$, where $d(v)$ denotes the number of adjacent vertices of $v$ in $V(G)$.
For $x, y \in V(G)$, distance $d(x, y)$ represents the length of the shortest path between $x$ and $y$, and $d(x, H)=min\{d(x,y)|y\in V(H)\}$ for $H\leq G$.
We denote by $K_{m,n}$ the complete bipartite graph.
Let $A(G)$ denote the adjacency matrix of graph $G$, which is a square matrix and $a_{ij}=1$ if and only if $v_i\sim v_j$, otherwise $a_{ij}=0$.

Let $m_G(\mu)$ denote the multiplicity of the eigenvalue $\mu$ in the adjacency matrix $A(G)$.
In 2013, Wong \cite{Wong1} characterize graphs with maximum diameter among all connected graph with rank $n$.
In 2022, Wang \cite{Wang1} characterized the graphs for which $m_G(0) = n - d$ using the method of star complements. 
In the same year, Chang and Li \cite{Chang1} characterized the graphs for which $m_G(0) = n - g - 1$.
In 2023, Du and Fonseca \cite{Du1} extend the definition of minimal graphs on
adjacency matrices to real symmetric matrices. They characterize all the trees for which
there is a real symmetric matrix with nullity $n-d$ and $n-d-1$.
In 2024, Zhang and Wong \cite{ZY} provided a complete characterization of the graphs for which $m_G(-1) = n - d$.
Therefor, we investigate the situation $m_G(-1)=n-d-1$.

Let $G$ be a graph and $u, v \in V (G)$.
We introduce (see \cite{MP}) a binary relation $\rho$ in $V(G)$ in the following way: $u\rho v$ if and only if $N_G[u] = N_G[v]$.
The relation $\rho$ is symmetric and transitive.
By this relation the vertex set $V(G)$ is divided to $k$ disjoint subsets $C_1$, $C_2$, $\ldots$, $C_k$, ($1 \leq k \leq n$, $n =|V (G)|$), such that every graph induced by the set $C_i$ ($i = 1, \ldots, k$) is a complete graph.
The graph $G^{c}$, obtained from $G$ by identification of all vertices from the same subset $C_i$ ($i = 1, \ldots, k$) is called the “$C$-canonical graph” of the graph $G$.
If $V(G)$ does not have any pair of vertices which lie in relation $\rho$, then $G^{c} = G$ and we say that $G$ is a $C$-canonical graph. For instance, all connected bipartite graphs distinct of $P_2$ are $C$-canonical graphs..

\section{Preliminaries}
\begin{Lemma}\cite{AB}
Let $v$ be a vertex of $G$, then $m_G(\mu)-1\leq m_G(\mu)\leq m_G(\mu)+1$
\end{Lemma}

The eigenvalues of $P_n$ are $\{2cos\frac{i\pi}{n+1}|i = 1, 2, ..., n\}$, we have the following lemma.

\begin{Lemma}
$m_{P_{n+1}}(-1)\leq1$, further, $m_{P_{n+1}}(-1)=1$ if and only if $n\equiv1 \ (\text{mod}\ 3)$
\end{Lemma}

\begin{Lemma}(\cite{PR}, Lemma 2.3)
Let $X$ be a star set for $\mu$ in $G$ and let $U$ be a proper subset of $X$.
Then $X-U$ is a star set for $\mu$ in $G-U$.
\end{Lemma}

\begin{Lemma}(\cite{DC}, Proposition 5.1.4)
Let $X$ be a star set for $\mu$ in $G$, and let $X =V (G)\backslash X$.
If $\mu=0$, then $X$ is a dominating set for $G$.
\end{Lemma}

\begin{Lemma}(\cite[3.1]{ZY})
Let $H$ be an induced subgraph of a connected graph $G$ for which $rank(A(H)+
I)\geq rank(A(G)+I)-1$, and $v\in V (G)\backslash V (H)$.
If $v\sim h$ and $N_H(v) = N_H[h]$ for a vertex $h\in V (H)$, then $N_G[v] = N_G[h]$.
\end{Lemma}

\begin{Lemma}(\cite[3.2]{ZY})
Let $H$ be an induced subgraph of a connected graph $G$ for which $rank(A(H)+
I)\geq rank(A(G)+I)-1$, and $u, v\in V (G)\backslash V (H)$.
If $u\sim v$ and $N_H(v) = N_H(h)$, then $N_G[v] = N_G[h]$.
\end{Lemma}

\begin{Lemma}(\cite[Theorem 5]{MP})
Let $G$ be a graph.
If two vertices $u$, $v$ satisfy $u\rho v$, then $m_G(-1)=m_{G-u}(-1)+1=m_{G-v}(-1)+1$.
\end{Lemma}

\begin{Lemma}(\cite[3.4]{ZY})
Let $G^{c}$ be the $C$-canonical graph of a graph $G$.
If $G$ is not a complete graph, then $d(G^{c})\geq1$ and $d(G^{c})=d(G)$.
\end{Lemma}

\begin{Lemma}
Let $G^{c}$ be the $C$-canonical graph of a graph $G$.
If $G$ is not a complete graph, then $m_G(-1)=n(G)-d(G)-1$ if and only if $m_{G^{c}}(-1)=n(G^{c})-d(G^{c})-1$.
\end{Lemma}

\begin{proof}
If $G$ is not a complete graph, by Lemma 2.8, $d(G^{c})=d(G)$.
Meanwhile, by Lemma 2.7, $m_G(-1)=m_{G^{c}}(-1)+n(G)-n(G^{c})$.
So, $m_G(-1)=n(G)-d(G)-1$ if and only if $m_{G^{c}}(-1)+(n(G)-n(G^{c})=n(G)-d(G)-1=n(G)-d(G^{c})-1)$.
Then we have $m_G(-1)=n(G)-d(G)-1$ if and only if $m_{G^{c}}(-1)=n(G^{c})-d(G^{c})-1$.
\end{proof}

By Lemmas 2.8 and 2.9, we know that for the eigenvalue $-1$, $G$ and $G^{c}$ possess identical properties; therefore, in this article, we will only consider the case where $G = G^c$.

\begin{Lemma}\cite[2.2]{WXL}
Let $G$ be a graph obtained from a given graph $H$ and $P_3$, disjoint with $H$, by identifying one pendant vertex of $P_3$ with a vertex $w$ of $H$. Then $m_G(-1)=m_{G-w}(-1)$
\end{Lemma}

\begin{Lemma}
If  $m_G(-1) = n - d - 1$ and $P_{d+1}$ is the diameter path of $G$ with $P_{d+1} \leq H \leq G$, then it follows that $rank(A(H) + I) \geq rank(A(G) + I) - 1$.
\end{Lemma}
\begin{proof}
If -1 is not an eigenvalue of $P_{d+1}$, then by the interlacing theorem, we have $rank(H+I)= rank(P_{d+1})$.
If -1 is an eigenvalue of $P_{d+1}$, then $rank(H+I) = rank(P_{d+1} + I)$ or $rank(P_{d+1} + I)+1$.
This completes the proof.
\end{proof}

For $x\in R^{n}$, we denote by $x_u$ the component of $x$ corresponding to $u\in V(G)$. 
With respect to a subset $U$ of $V(G)$, we define 

\centerline{ $\mathbb{Z}(U)=\{x\in R^{n}: x_u=0, \forall u\in U\}$.} 

\begin{Lemma}\cite[Lemma 2.1]{WXL}
Let $G$ be a graph of order $n$ and $U$ be a subset of $G$.
If $V(G,\mu))\cap Z(U)=0$, then $m_G(\mu)\leq |U|$.
\end{Lemma}

\begin{Lemma}\cite{PR1}
Let $G$ be a graph, $V(G)=\{1, 2, ..., n\}$, $X\subseteq V (G)$. 
Let $F$ be the matrix which represents the orthogonal projection of $R^{n}$ onto the eigenspace $\mathcal{E}_{A(G)}(\mu)$ with respect to the standard orthonormal basis $\{e_1, e_2,\ldots, e_n\}$ of $R^{n}$.
Then

(i) $X$ is a star set for $\mu$ in $G$ if and only if the vectors $Fe_i$, $i\in X$, form a basis for $\mathcal{E}_{A(G)}(\mu)$.

(ii)The matrix $F$ is a polynomial in $A$ and $\mu Fe_v = \Sigma_{u\sim v}Fe_u$.
\end{Lemma}

\section{$m_P(-1)=0$ and $m_G(-1)=n-d-1$}
\begin{figure}[H]
  \centering
  \includegraphics[width=1.0\linewidth]{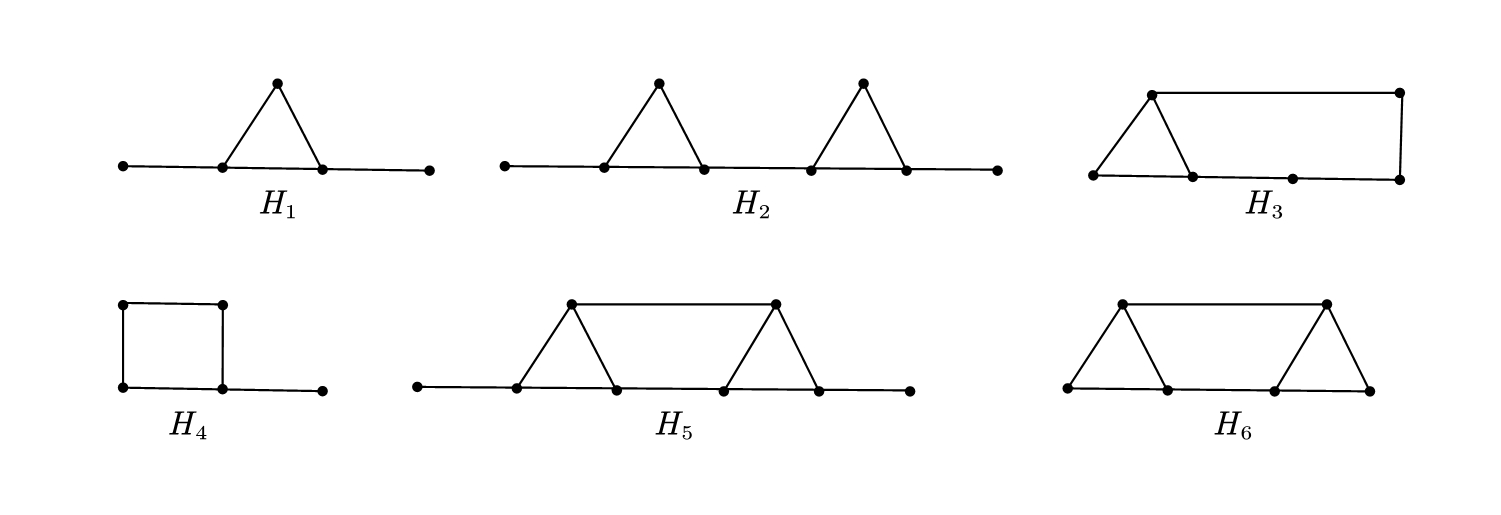}

\end{figure}

After performing computations using a computer, it can be concluded that $m_{H_i}(-1)=0$ for $1\leq i\leq6$.

\begin{figure}[H]
  \centering
  \includegraphics[width=1.0\linewidth]{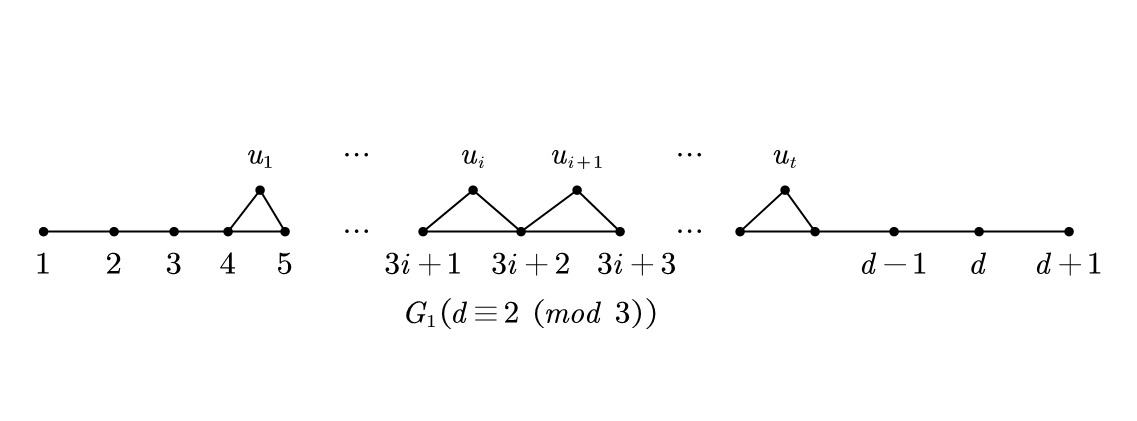}

\end{figure}

\begin{figure}[H]
  \centering
  \includegraphics[width=1.0\linewidth]{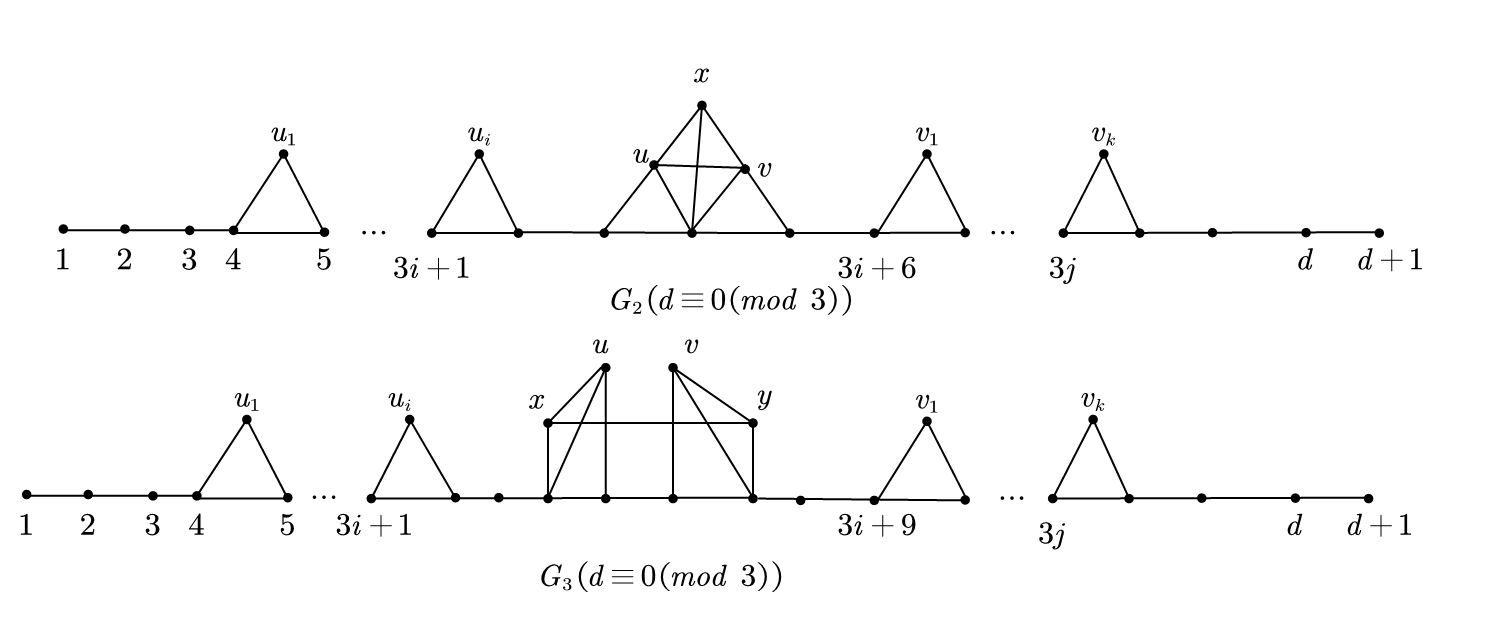}

\end{figure}

\begin{Theorem}
Let $G$ be a connected graph and $P = P_{d+1}=v_1\sim v_2\sim\cdots \sim v_{d+1}$ be the diameter path of $G$.
Given that -1 is not an eigenvalue of $P$, we have $m_G(-1) \leq n - d - 1$.
The equality holds if and only if $P \leq G \leq G_i $ for $1 \leq i \leq 3$.
\end{Theorem}

\begin{proof}
Since -1 is not an eigenvalue of $P$, we can deduce from Lemma 2.1 that $m_G(-1) \leq n - d - 1 $, with equality holding if and only if $X = V(G - P)$ forms a star set.
From Lemma 2.4, we know that $P$ is a dominating set for $G$, which implies $N_G(P) = V(G - P)$.
Let $N_G(P) = \cup_{i\geq1} N_i$ and $N_i=\{x: |N_P(x)| = i\}$.

\textbf{Claim 1}: $N_G(P)=N_1\cup N_2$

If $i > 3$ and $x \in N_i $, where $x$ is adjacent to $v_{j_1}$, $v_{j_2}$, $v_{j_3}$, $v_{j_4}$ with $j_1 < j_2 < j_3 < j_4$, then the path $v_1 v_2 \ldots v_{j_1} x v_{j_4} \ldots v_{d+1}$ forms a path shorter than $P$, which contradicts the fact that $P$ is a diameter path.
If $i = 3$ and $x \in N_3$, with $x$ being adjacent to $v_{j_1}$, $v_{j_2}$, $v_{j_3}$, it follows that $j_3 = j_1 + 2 = j_2 + 1$.
If this condition is not satisfied, then $v_1 v_2 \ldots v_{j_1} x v_{j_3} \ldots v_{d+1}$ would again represent a path shorter than $P$, leading to a contradiction.
When $j_3 = j_1 + 2 = j_2 + 1$, we find that $N_P(x) = N_P[v_{j_2}]$.
According to Lemma 2.5, this implies  $N_G[x] = N_G[v_{j_2}]$, which contradicts the assumption that $G = G^{c}$.
Therefore, we conclude that $N_G(P) = N_1 \cup N_2$.

\textbf{Claim 2}: If $x,y\in V(G-P)$, then $N_P(x)\neq N_P(y)$.

If $x \sim y$, then $N_P[x] = N_P[y]$. 
According to Lemmas 2.6 and 2.11, it follows that $N_G[x] = N_G[y]$, which contradicts the condition $G = G^c$. 
If $x \nsim y$, consider $P + x + y$. 
By Lemma 2.13, we have $Fe_x = Fe_y$, which contradicts the fact that $Fe_x$ and $Fe_y$ are linearly independent. 
Therefore, we conclude that $N_P(x) \neq N_P(y)$.

Since $m_P(-1) = 0$, we can infer from Lemma 2.2 that $d\equiv 0$ or $2 \ (\text{mod} \ 3)$ . 
Therefore, we will discuss two cases in the following analysis.

\textbf{Case 1}: $d\equiv2 \ (\text{mod} \ 3)$

\textbf{Claim 1.1}: $N_1=\emptyset$.

First, by the interlacing theorem, we know that $m_{P+u}(-1) = 1$.

If $N_1 \neq \emptyset$, let $x \in N_1$ such that $x \sim v_m$. 
If $m \equiv 0$ or $1 \mod 3 $, then by Lemma 2.10, we have $m_{P+x}(-1) = 0$, leading to a contradiction. 
If $m \equiv 2 \mod 3$, it follows that $m_{P+x}(-1) = m_{K_{1,3}} = 0$, which also leads to a contradiction. 
Therefore, we conclude that $N_1 = \emptyset$.

\textbf{Claim 1.2}: If $u \in N_2$ and $u$ is adjacent to $v_i$ and $v_j$ with $i < j$, then it follows that $j = i + 1$.

If $j - i > 2$, then the path $v_1 v_2 \ldots v_i u v_j \ldots v_{d+1}$ is shorter than $P$, leading to a contradiction.

When $ j = i + 2$ and $i \equiv 0 \ (\text{mod} \ 3)$, then $m_{P+u}(-1) = 0$, resulting in another contradiction.

When $j = i + 2$ and $i \equiv 1 \ (\text{mod} \ 3)$, then $m_{P+u}(-1) = m_{H_4}(-1) = 0 $, resulting in yet another contradiction.

The case for $j = i + 2$ and $i \equiv 2 \ (\text{mod} \ 3)$ can be similarly analyzed, leading to a contradiction.
In summary, we conclude that $j = i + 1$.

\textbf{Claim 1.3}: If $u \in N_2$ and $u$ is adjacent to $v_m$ and $v_{m+1}$, then $m \equiv 0$  or  $1 \ (\text{mod} \ 3)$.

When $m \equiv 0 \ (\text{mod} \ 3)$, then by using Lemma 2.10, we obtain $m_{P+u}(-1) = 0 $, which leads to a contradiction.

When $m \equiv 1 \ (\text{mod} \ 3)$, then $m_{P+u}(-1) = 1 $, satisfying the condition.

When $m \equiv 2 \ (\text{mod} \ 3)$, then $m_{P+u}(-1) = 1$, satisfying the condition.
Thus, the proof is established.

\textbf{Claim 1.4}: If $u, v \in N_2$ with $u \sim v_m, v_{m+1}$ and $v \sim v_h, v_{h+1}$, we may assume $m < h$. 
Then, $u \not\sim v$ and $P\leq G\leq G_1$.

First, let us assume that $u \sim v$. 
When $m \equiv 1 \ (\text{mod} \ 3)$, since $P$ is a diameter path, we only need to consider the case where $h = m + 1$.
In this case, $m_{P + u + v}(-1) = 1 $, which leads to a contradiction.
Similarly, when $m \equiv 2 \ (\text{mod} \ 3)$, we only need to consider the case where $h = m + 2$. 
In this situation, $m_{P + u + v}(-1) = m_{H_5}(-1) = 0$, which also results in a contradiction.
Thus, it follows that $u \not\sim v$.

If $ m \equiv 1 \ (\text{mod} \ 3)$ and $h \equiv 2 \ (\text{mod} \ 3)$, then $m_{P + u + v}(-1) = 2 $, which satisfies the conditions.
If $m \equiv 2 \ (\text{mod} \ 3)$ and $h \equiv 1 \ (\text{mod} \ 3) $, then $m_{P + u + v}(-1) = m_{H_2}(-1) = 0 $, leading to a contradiction.
If both $m \equiv 1 \ (\text{mod} \ 3) $ and  $h \equiv 1 \ (\text{mod} \ 3) $ (or if $m \equiv 2 \ (\text{mod} \ 3) $ and $ h \equiv 2 \ (\text{mod} \ 3) )$, we can apply Lemmas 2.7 and 2.10, which yield $ m_{P + u + v}(-1) = 2 $, satisfying the conditions.
Therefore, based on the above discussion, $G_1$ is a maximal graph that satisfies the given conditions.
Furthermore, by reapplying Lemmas 2.7 and 2.10, we obtain $m_{G_1}(-1)=t= n - d - 1$, which concludes the proof.

We will now begin to consider the case where $d \equiv 0 \ (\text{mod} \ 3)$, while omitting certain proofs that are similar to those in Case 1.

\textbf{Caes2}: $d\equiv0 \ (\text{mod} \ 3)$

\textbf{Claim 2.1}: If $x\in N_1$ and $x\sin v_m$, then $m \equiv 1 \ (\text{mod} \ 3)$.

If $m \equiv 0$ or $2$ $\ (\text{mod} \ 3)$, then $m_{P+u}(-1)=0$, leading to a contradiction. 

If $m \equiv 1 \ (\text{mod} \ 3)$, then $m_{P+u}(-1)=1$, satisfying the condition.

\textbf{Claim 2.2}: If $x, y \in N_1$, then $x \sim y$. 
Therefore, $|N_1| \leq 2$.

Let $x \sim v_m$ and $y \sim v_h$.
According to Claim 2.1, we know that $m = 3t + 1$ and $h = 3l + 1$.
If $x \not\sim y$, then $m_{P + x + y}(-1) = m_{P_3}(-1) = 0$, which leads to a contradiction.
If $x \sim y$, we may assume $m < h$, thus $h = m + 3$.
Otherwise, the sequence $v_1 v_2 \ldots v_m xy v_h \ldots v_{d+1}$ would constitute a path shorter than $P$.
In this case, we have $m_{P + x + y}(-1) = m_{C_6}(-1) = 2$, which satisfies the conditions.

\textbf{Claim 2.3}: If $u \in N_2$ and $u$ is adjacent to $v_i$ and $v_j$ with $i < j$, then it follows that $j = i + 1$.

\textbf{Claim 2.4}: If $u \in N_2$ and $u$ is adjacent to $v_m$ and $v_{m+1}$, $m \equiv 0$  or  $1 \ (\text{mod} \ 3)$.

When $m \equiv 0 \ (\text{mod} \ 3)$, then by using Lemma 2.10, we obtain $m_{P+u}(-1) = m_{P_2}(-1) = 1 $, satisfying the condition.

When $m \equiv 1 \ (\text{mod} \ 3)$, then  $m_{P+u}(-1) = m_{P_2}(-1) = 1$, satisfying the condition. 

When $m \equiv 2 \ (\text{mod} \ 3)$, then $m_{P+u}(-1) = m_{H_1} = 0$, leading to a contradiction;

\textbf{Claim 2.5}: Let $u, v \in N_2$ such that $u \sim v_{3a}, v_{3a + 1}$, and $v \sim v_{3b + 1}, v_{3b + 2}$.
If $b \geq a$, then it follows that $u \sim v$ and $b = a$.
Conversely, if $b < a$, then $u \nsim v$.

If $b \geq a$ and $u \not\sim v$, then according to Lemma 2.10, we have $m_{P + u + v}(-1) = 0 < 2$, which leads to a contradiction.
Furthermore, since $P$ is a diameter path, we can conclude that $a = b$.

If $b < a$ and $u \sim v$, considering that $P$ is a diameter path, we have $b = a - 1$.
Thus, $m_{P + u + v}(-1) = m_{H_6}(-1) = 0 < 2 $, which leads to a contradiction.
Therefore, it follows that $u \nsim v$.

\textbf{Claim 2.6}: If $x \in N_1$, $u \in N_2$, $x \sim v_i$, $u \sim v_j, v_{j+1}$, and $u \not\sim x $, then if $ j > i $, it follows that $j \equiv 0 \ (\text{mod} \ 3) $.
Conversely, by symmetry, if $j < i$, then $j \equiv 1 \ (\text{mod} \ 3) $.

From Claim 2.4, we know that $j \equiv 0$ or $1 \ (\text{mod} \ 3)$.
If $j > i$ and $j \equiv 1 \ (\text{mod} \ 3)$, then $m_{P + x + u}(-1) = m_{P_3}(-1) = 0$, which leads to a contradiction.
On the other hand, if $j > i$ and $j \equiv 0 \ (\text{mod} \ 3) $, then $m_{P + x + u}(-1) = m_{C_3}(-1) = 2$, satisfying the conditions.
By symmetry, the proof is completed.

\textbf{Claim 2.7}: Let $x \in N_1$ and $u \in N_2 $, with $x \sim v_i$, $u \sim v_j$, $ v_{j+1} $.
If $x \sim u$, then it follows that $j = i$ or $j = i - 1$.

By symmetry, it is sufficient to prove that when $j \geq i$, then $j = i$.
If $j \equiv 1 \ (\text{mod} \ 3)$, then $i = j$; otherwise, $v_1 v_2 \ldots v_i x u v_{j+1} \ldots v_{d+1}$ would form a path shorter than $P$ ; when $i = j$, we have $m_{P + u + v}(-1) = m_{C_3}(-1) = 2$, which satisfies the conditions.
If $j \equiv 0 \ (\text{mod} \ 3)$, then $j = i + 2$; otherwise, a path shorter than $P$ would emerge.
When $j = i + 2$, we find $m_{P + u + v}(-1) = m_{H_3}(-1) = 0$, leading to a contradiction.
Therefore, the conclusion is established.

\textbf{Claim 2.8}: Let $N_1 = \{ x, y \}$.
If there are two vertices  $u, v\in N_2$ that are adjacent to vertices in $N_1$, then $P + u + x + y\cong H_7$.

According to claims 2.6 and 2.7, when there are two vertices $u, v\in N_2$ that are adjacent to vertices in $N_1$, it follows that $P + u + v + x + y \cong H_7 $.
Based on Lemmas 2.7 and 2.10, we have $m_{H_7}(-1) = 2 + m_{C_6}(-1) = 4$, which satisfies the conditions and thereby concludes the proof.
\begin{figure}[H]
  \centering
  \includegraphics[width=1.0\linewidth]{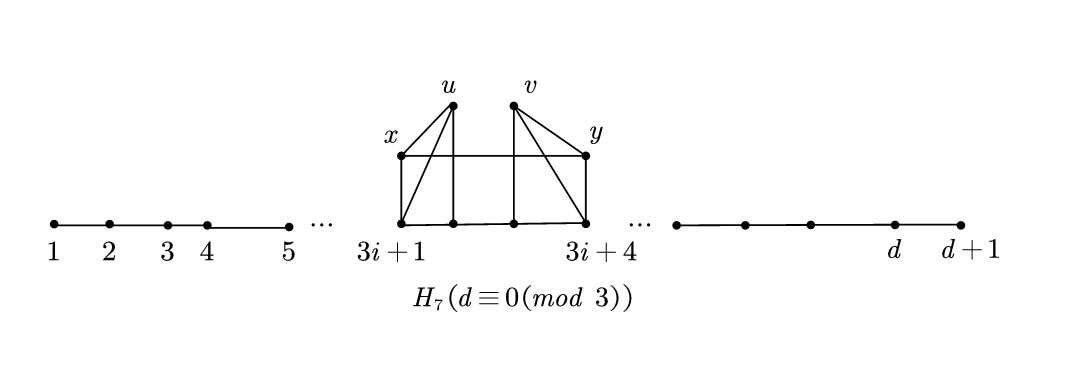}

\end{figure}

\textbf{Claim 2.9}:  When $d \equiv 0 \ (\text{mod} \ 3) $, $P\leq G\leq G_2$ or $G_3$.

Based on the above discussion, it can be concluded that $G_2$ or $G_3$ are maximal graphs that satisfy the given conditions.
Furthermore, according to Lemmas 2.7 and 2.10, we obtain $m_{G_2}(-1) = m_{G_3}(-1) = i + k + 4 = n - d - 1$, which satisfies the conditions. The proof is thereby established by the alternating theorem.
\end{proof}
\section{$m_P(-1)=1$ and $m_G(-1)=n-d-1$ }

\begin{Theorem}
Let $G$ be a connected graph and $P = P_{d+1}=v_1\sim v_2\sim\cdots \sim v_{d+1}$ be the diameter path of $G$.
If $m_P(-1)=1$ and $m_G(-1)=n-d-1$, then $P<G\leq G_4$ and $N_1\neq \emptyset$.
\end{Theorem}

\begin{figure}[H]
  \centering
  \includegraphics[width=1.0\linewidth]{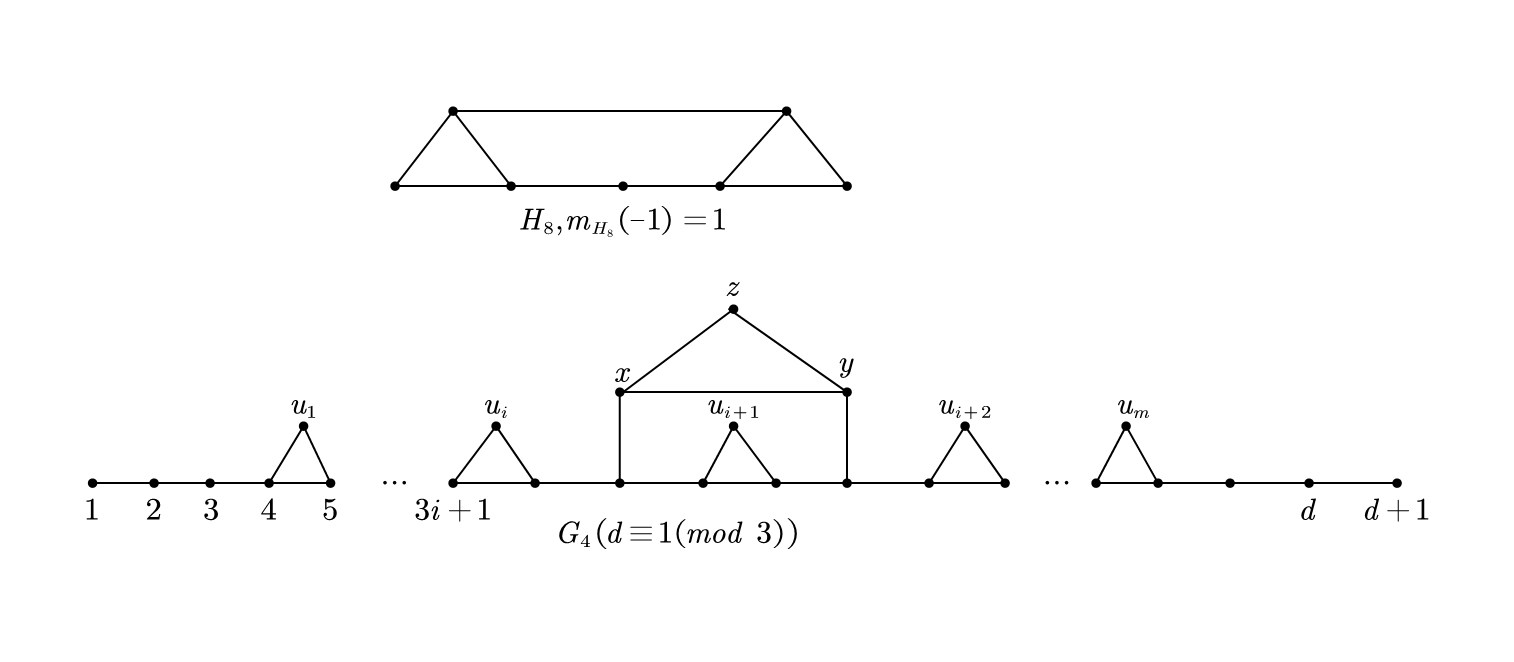}

\end{figure}

\begin{proof}

Initially, since $m_P(-1) = 1$, we conclude from Lemma 2.2 that $d = 3k + 1$.
Furthermore, there cannot exist a point $z$ that is at a distance of 3 from $P$. To illustrate this, suppose ( z ) is a point at a distance of 3 from $P$, $y$ is a vertex at a distance of 2 from $P$, and $x$ is adjacent to  $P$.
According to Lemma 2.10, we then have $m_{P + x + y + z}(-1) = m_P(-1) = 1 < 2$, which leads to a contradiction.

\textbf{Claim 1}: If $x \in N_1$ and $x \sim v_i$, then $m_{P + x}(-1) = 1$ and $i \equiv 0 \mod{3}$.

According to Lemma 2.11, we have $m_{P+x}(-1) \leq 1 + \dim(V(G, -1) \cap Z(v_1) $.
Since $\dim(V(G, -1) \cap Z(v_1)) = 0$, it follows that $m_{P+x}(-1) \leq 1$.
If $m_{P+u}(-1) = 0$, then it would imply that $m_G(-1) \leq n - d - 2$, which leads to a contradiction.
Therefore, we conclude that $m_{P+u}(-1) = 1$.
This also indicates that $V(G - P - x)$ is a star set.

If $i \equiv 0 \ (\text{mod}\ 3)$, then by utilizing Lemma 2.10, we obtain $m_{P+x}(-1) = m_{P_2}(-1) = 1$, which satisfies the condition.
If $m \equiv 1$  or  $2\ (\text{mod}\ 3) $, then $m_{P+u}(-1) = 0$, leading to a contradiction.

\textbf{Claim 2}: If $x, y \in N_1$ and $x \sim v_{3a}$, $y \sim v_{3b}$, then $a \neq b$ and $u \sim v $.
Hence, it follows that $|N_1| \leq 2$.

First, we demonstrate that $a \neq b$.
If $a = b$ and $x \sim y$, then $N_P(x) = N_P(y)$.
By Lemma 2.6, this implies $N_G[x] = N_G[y]$, which leads to a contradiction.
Alternatively, if $a = b$ and $x \not\sim y$, then $m_{P+x+y}(-1) = 1 < 2$, resulting in another contradiction.

If $x \not\sim y$, then by Claim 1 and Lemma 2.10, we have $m_{P+x+y} = m_{P_2}(-1) = 1<2$, which leads to a contradiction.
Hence, it follows that $x \sim y$.
Furthermore, considering that $P$ is a diameter path, we conclude that $|N_1| \leq 2$.

\textbf{Claim 3}: Let $u \in N_2$ and $u \sim v_m, v_{m+1}$.
Then $m \equiv 1 \ (\text{mod} \ 3)$ and $m_{P+u}(-1) = 2$.

First, we know that $1 \leq m_{P+u}(-1) \leq 2$.
If $m \equiv 0 \ (\text{mod} \ 3)$  or  $m \equiv 2 \ (\text{mod} \ 3)$ , then $ m_{P+u}(-1) = 0 $, which leads to a contradiction.
If $m \equiv 1 \ (\text{mod} \ 3) $, then $m_{P+u}(-1) = m_{C_3}(-1) = 2 $, thus proving the claim.
Similar to Claim 2, for $u, v \in N_2$, we have $N_P(u)\neq N_P(v)$.

\textbf{Claim 4}: Let $x \in N_1$ and $u, v \in N_2$. Then $x \not\sim N_2$ and $u \not\sim v$.

Assume $u \sim v$.
According to Claim 3, we may assume $u \sim v_{3a+1}$, $v_{3a+2}$, and $v \sim v_{3b+1}$, $v_{3b+2}$.
Considering that $P$ is a diameter path, we have $|b-a|=1$.
Thus, $m_{P+u+v}(-1) = m_{H_8}(-1) = 1 < 2 $, which leads to a contradiction.
If $x \sim N_2$, we may assume $x \sim v_{3c}$ and $u$.
Given that $P$ is a diameter path, we have $|b-c|=1$.
Therefore, $m_{P+x+u}(-1) = 1 < 2$, which results in another contradiction.

Next, we begin to consider the vertices that are at a distance of 2 from $P$.

\textbf{Claim 5}: If $d(z, P) = 2$, then $z \sim N_1$ and $z \not\sim N_2$.

If $x \in N_1$ and $z \sim x $, then $m_{P+x+z}(-1) = 2$, satisfying the condition.
On the other hand, if $u \in N_2$ and $z \sim u$, then $m_{P+u+z}(-1) = 1 < 2 $, which leads to a contradiction.

\textbf{Claim 6}: Let $N^{2} = \{ z \mid d(z, P) = 2 \}$.
Then it follows that $|N^{2}| \leq 1 $.

If $z, w \in N$, by Claim 5, we know there exist $x, y \in N_1$ such that $x \sim z$ and $y \sim w$.
If $x = y$ and $z \sim w$, let $H = P + x$; then $N_H(z) = N_H(w)$.
By Lemma 2.5, we have $N_G[z] = N_G[w]$, which leads to a contradiction.
If $x = y$ and $z \not\sim w$, then $m_{P+x+z+w}(-1) = 1 < 3$, which is also a contradiction.
Thus, we conclude that $x \neq y$.
At this point, based on Claim 1 and Claim 2, we can assume without loss of generality that $x \sim v_{3a}$ and $y \sim v_{3a+3}$.
If $z \not\sim y$, then $m_{P+x+y+z}(-1) = 1 < 2$, which is again a contradiction.
If $z \sim y$, then $m_{P{x+y+z}}(-1) = 1 + m_{C_3}(-1) = 3$, satisfying the condition, which implies $z \sim x, y $.
Similarly, we can conclude that $w \sim x, y$.
If $z \sim w$, let $H = P + x + y )$; then $N_H(z) = N_H(w)$.
By Lemma 2.6, we know $N_G[z] = N_G[w]$, leading to a contradiction.
If $z \not\sim w$, then $m_{P+x+y+z+w}(-1) = 2 < 4$, which is also a contradiction.
In summary, we have $|N^2| \leq 1 $. If $z \in N^2$ and $ x, y \in N_1 $, then it follows that $z \sim x, y $.

\textbf{Claim 7}: $P<G\leq G_4$ and $N_1\neq \emptyset$.

First, based on the discussion above, we know that $G_4$ is a maximal graph that satisfies the conditions.
Furthermore, by applying Lemmas 2.7 and 2.10, we can conclude that $m_{G_4}(-1) = m + 3 = n - d - 1$, which satisfies the required conditions. If $N_1 = \emptyset$, then $V(G) = V(P) \cup N_2$, leading to $m_G(-1) = n - d $, which presents a contradiction.
Therefore, we conclude that $N_1 \neq \emptyset$, and the theorem is proven.

\end{proof}

\end{document}